\documentclass[amstex,12pt, amssymb]{article}

\usepackage{mathtext}
\usepackage[cp1251]{inputenc}
\usepackage[T2A]{fontenc}
\usepackage[dvips]{graphicx}
\usepackage{amsmath}
\usepackage{amssymb}
\usepackage{amsxtra}
\usepackage{latexsym}
\usepackage{ifthen}

\newcounter{lemma}[section]

\newcounter{corollary}[section]

\newcounter{remark}[section]

\newcounter{theorem}[section]

\newcounter{proposition}[section]

\numberwithin{equation}{section}

\def\XXint#1#2#3{{\setbox0=\hbox{$#1{#2#3}{\int}$}
     \vcenter{\hbox{$#2#3$}}\kern-.5\wd0}}

\def\cc{\setcounter{equation}{0}
\setcounter{figure}{0}\setcounter{table}{0}}

\begin{document}

\markboth{\centerline{VLADIMIR RYAZANOV}} {\centerline{PROOF OF RIEMANN HYPOTHESIS}}

\author{{Vladimir Ryazanov}}

\title{{\bf  A proof of the Riemann hypothesis\\
on zeros of $\zeta -$function.}}

\maketitle

\large \begin{abstract} In his famous presentation at the International Congress of
Mathematicians held in Paris in 1900, David Hilbert included the
Riemann Hypothesis on zeros of $\zeta
-$function as number 8 in his list of 23 challenging
problems published later. After over 150 years, it is one of the few
on that list that have not been solved. At present many
ma\-the\-ma\-ti\-cians consider it the most important unsolved
problem in mathematics.

Recall that, exactly one hundred years later, the Clay Mathematics
Institute has published a list of 7 unsolved problems for the 21st
century, including 6 unresolved problems from the Hilbert list,
offering a reward of one million dollars for a solution to any of
these problems.

One of them is the {\bf Riemann hypothesis}, i.e., a conjecture that
the so-called Riemann zeta function has as its zeros only complex
numbers with real part $1/2$ in addition to its trivial zeros at the
negative even integers. It was proposed by Bernhard Riemann in his
1859 paper \cite{R}. The Riemann zeta function plays a great role in
analytic number theory as well as in physics, probability theory and
quantum mechanics.

In this preprint, applying the known Beurling--Nyman
criterion, it is just proved the Riemann hypothesis on zeros of $\zeta -$function.
\end{abstract}


\bigskip

{\bf Dedicated to my wife Tamara}

\section{Introduction}

The {\bf Riemann zeta function} $\zeta(s)$ is a function of a
complex variable $s$ that analytically continues the sum of the {\bf
Dirichlet series}
\begin{equation}\label{eqDEF}
\zeta(s)\ =\ \sum\limits_{n=1}\limits^{\infty}\ \frac{1}{n^s}\ , \ \
\ \ \ \  \mbox{Re}\ s\ >\ 1\ .
\end{equation}
As known, the series (\ref{eqDEF}) is extended to the meromorphic
function $\zeta(s)$ of the whole plane having only one simple pole
at the point $s=1$.

\medskip

The Riemann zeta function satisfies the {\bf Riemann functional
equation}
\begin{equation}\label{eq}\zeta (s)\ =\ 2^{s}\pi ^{s-1}\ \sin \left({\frac {\pi s}{2}}\right)\
\Gamma (1-s)\ \zeta (1-s)\ \ \ \ \ \  \forall\ s\in\mathbb C
\end{equation}
which is an equality of meromorphic functions where $\Gamma(s)$ is
the {\bf gamma function of Euler}, see \cite{R}, see also \cite{T}.
Recall that $\Gamma(s)$ is a meromorphic function on the whole
complex plane $\mathbb C$ having no zeros and only simple poles at
the points $s=0,-1,-2,\ldots $. Equation (\ref{eq}) implies that
$\zeta(s)$ has simple zeros at all even negative integers $s = -2n$,
these are the {\bf trivial zeros of} $\zeta(s)$.

\medskip

Riemann has also found in \cite{R} a symmetric form of the
functional equation (\ref{eq}). One of such equivalent forms, see
e.g. \cite{L} and also \cite{T}, is the equation
\begin{equation}\label{eqSYM}
\xi(s)\ =\ \xi(1-s)\ \ \ \ \ \  \forall\ s\in\mathbb C
\end{equation}
where
\begin{equation}\label{eqXI}\xi (s)\ =\ \frac{1}{2}\ \pi
^{-\frac{s}{2}}\ s(s-1)\ \Gamma \left(\frac{s}{2}\right)\ \zeta (s)\
.
\end{equation}
Note that by the previous items the function $\xi(s)$ is an entire
function, i.e., an analytic function in the whole complex plane
$\mathbb C$ without any poles, and, moreover,  $\xi(s)$ has no above
trivial zeros of $\zeta(s)$ but all their rest zeros coincide.

\bigskip

{\bf Remark 1.} After the replacement $z=s-1/2$, equation
(\ref{eqSYM}) can be written in the form
\begin{equation}\label{eqHALF}
\xi\left(\frac{1}{2} + z\right)\ =\ \xi \left(\frac{1}{2} -
z\right)\ \ \ \ \ \ \forall\ z\in\mathbb C
\end{equation}
meaning that the function $\xi(s)$ is symmetric with respect to the
point $z_0 = 1/2$. Thus, to verify the Riemann hypothesis it is
sufficient to prove the absence of zeros of $\zeta(s)$ in the
half-plane $\rm{Re}\, s>1/2$. Note also by the way that zeros of the
Riemann zeta function are symmetric with respect to real axes
because at all $\overline{\zeta(s)}=\zeta(\bar s)$.


Note that there exist fine monographs devoted to the
theory of the Riemann zeta function, see e.g. \cite{Ap}, \cite{Ay}, \cite{E}, \cite{HW},
\cite{Iv}--\cite{KV}, \cite{Lap} and the classic \cite{T}. Moreover,
it was even appeared the 3 volumes of equivalents of the Riemann
hypothesis, see \cite{B1}--\cite{B3}. The great number of such
equivalents makes possible, on the one hand, to attack the Riemann
hypothesis from many positions and, on the other hand, to obtain
many consequences in the case of its proof or its disproof. We
prefer one of these equivalents.


\cc
\section{The Beurling--Nyman criteria}

Let us recall the contents of the paper \cite{Beu} of the known
Swedish mathematician Arne Beurling. Denote by $\{ \tau\}$ the
fractional part $\tau - [\tau]$ of a real number $\tau$ where
$[\tau]$ is the greatest integer that is less or equal to $\tau$.
Denote also by $\frak{B}$ the collection of all functions $\varphi
:(0,1]\to\mathbb R$ of the form
\begin{equation}\label{eqB}
\varphi(t)\ =\ \sum\limits_{k=1}\limits^{n
}c_k\left\{\frac{\theta_k}{t}\right\}, \ \  c_k\in\mathbb R,\
\theta_k\in(0,1],\ k=1,\ldots , n,
\end{equation}
with the condition
\begin{equation}\label{zero}
\sum\limits_{k=1}\limits^{n }c_k\theta_k\ =\ 0\ .
\end{equation}
Now, let $\frak{B}_p$ be the closure of $\frak{B}$ in
$L_p=L_p(0,1)$, $1 < p <\infty$. It is shown in \cite{Beu} that
$\frak{B}_p = L_p$ if and only if the function $f(t)\equiv 1$,
$t\in(0,1)$, is in $\frak{B}_p$. Moreover, it is shown in \cite{Beu}
that the Riemann zeta function has no zeros in $\rm{Re}\, s > 1/p$,
$p\in(1,2]$, if and only if $\frak{B}_p = L_p$. Thus, by Remark 1 we
have from this the following consequences.

\medskip

{\bf Theorem A.} {\it The Riemann hypothesis is true if and only if
the function $f(t)\equiv 1$, $t\in(0,1]$, can be approximated in
$L_2$ by a sequence in the class $\frak{B}$.}

\medskip

{\bf Theorem B.} {\it The Riemann zeta function has no zeros in the
half--plane $\rm{Re}\, s > 1/p$ for $p\in(1,2]$ if and only if the
function $f(t)\equiv 1$, $t\in(0,1]$, can be approximated in $L_p$
by a sequence in the class $\frak{B}$.}

\medskip

Theorem A was first proved in the thesis \cite{Ny} of Bertil Nyman
(1950). Recall that Beurling was his supervisor. The paper
\cite{Beu} first (1955) represented and generalized this result as
Theorem B. Later on, the Beurling--Nyman criterion was reproved and
generalized in many different ways, as well as, the approach was
applied for the research of the problem on the distribution of zeros
of the Riemann zeta function, see e.g. monographs \cite{B2}, \cite{Iv}, \cite{Iw} and articles \cite{Baez1}--\cite{Ber},
\cite{Bur}--\cite{Ch}, \cite{H}, \cite{LR}--\cite{Nik}, \cite{Ro}, \cite{Ya2}--\cite{Ya3}.
It is impossible to list here thousands of other articles devoted immediately or indirectly to the Riemann hypothesis.
\medskip

Baez-Duarte in \cite{Baez4} proved Theorem 1 on a new version of
criteria A in terms of class $\frak {D}$ of functions $\varphi$ in
$L_2(0,\infty)$ of the form (\ref{eqB}) with the special
$\theta_k=1/k$, $k=1,2,\ldots$, generally speaking
without (\ref{zero}), see also \cite{Bag} and \cite{B2}.
\medskip

{\bf Theorem C.} {\it The Riemann hypothesis is true if and only if
the characteristic function of interval $(0,1]$ can be approximated
in $L_2$ by a sequence in class $\frak{D}$.}

\medskip

Theorem C admits wording more convenient for our purposes. Namely,
let us denote by $\frak{R}$ the collection of all functions $\varphi
:(0,\infty)\to\mathbb R$ of the form
\begin{equation}\label{eqBSTAR}
\varphi(t)\ =\ \sum\limits_{k=1}\limits^{n} h_k\, a_k(t)\ , \ \
h_k\in\mathbb R,\ \ a_k(t)\ :=\ \left\{\frac{t}{k}\right\}\ ,\
\hbox{$k$ and $n\in\mathbb N$}\ ,
\end{equation}
by $L_2^*$ the Hilbert space $L_2((0,\infty);\,t^{-2}dt)$, by $\|
f\|_*$ the norm of $f$ in $L_2^*$ and by $\langle f,g\rangle_*$ the
scalar product of functions $f$ and $g$ in $L_2^*\ $,
\begin{equation}\label{scalar}
\langle f,g\rangle_*\ :=\ \int\limits_0^{\infty}f(t)\,g(t)\
\frac{dt}{t^2}\ .
\end{equation}

Applying the replacement $t \mapsto 1/t$, we come to the following
equivalent formulation of the Baez-Duarte criterion that makes our
study more visual.

\medskip

{\bf Theorem D.} {\it The Riemann hypothesis is true if and only if
the characteristic function $\chi$ of $[1,\infty)$ can be
approximated by a sequence of the class $\frak{R}$ in $L^*_2$.}

\medskip

{\bf Remark 2.} Note that the functions $\varphi\in\frak{R}$ are
linear with the slope $\sum\frac{h_k}{k}$ on all open intervals of
length $1$ that appeared in $(0,\infty)$ after removing all natural
numbers, continuous
from the right at each point $k\in\mathbb N$, $\varphi(+0)=0$, and have the jump $-h_k$
at each point $k m$, $k,m\in\mathbb N$, where the corresponding
jumps are summarized if $k_1m_1=k_2m_2$ for some indexes $m_1\ne m_2$
and $k_1\ne k_2$.

\medskip

Let us\emph{} add a little bit algebra to this geometric interpretation of the problem.

\cc
\section{On the Connection to the M\"obius Function}

Recall that the {\bf M\"obius function} $\mu(n)$ of the natural parameter
$n=1,2,\ldots$ is defined as follows: (i) $\mu(1)=1$, (ii)
$\mu(n)=0$ if $n$ is divisible by a square of a prime $p>1$, (iii)
$\mu(n)=(-1)^k$ if $n$ is the product of $k$ distinct primes $p>1$.

\medskip

The function $\mu(n)$ occurs implicitly in the work of Euler as
early as 1748, but M\"obius, in 1832, was the first to begin systematically
studying its properties, see e.g. \cite{L}, p. 567--587 and
901.

\medskip

Its {\bf characteristic property}, see e.g. Theorem 263 in
\cite{HW}, see also Theorem 2.1 in \cite{Ap}, p. 24, Theorem 7.2
in \cite{Ay}, p. 103, and \cite{L}, p. 575, is the following:
\begin{equation}\label{eq2}\boxed{\ \sum\limits_{d\,|\,n}\mu(d)\ =\ 0\ }\ \ \ \forall\ n\ >\ 1 \ .\end{equation}

{\bf Remark 3.} Arguing by induction, it is easy to see that if some
arithmetical function $\alpha(n)$, $n\in\mathbb N$, satisfies this
property and $\alpha(1)=1$, then $\alpha(n)\equiv\mu(n)$, because we
have by (\ref{eq2}) that
\begin{equation}\label{eq3}{\ \alpha(n+1)\ =\ -\sum\limits_{d|n+1,\,d<n+1}\alpha(d)\ }\ \ \ \ \ \ \forall\ n\, \in\, \mathbb N \ ,\end{equation}
i.e., we have a guarantee for the next inductive step, and,
consequently, such a function is uniquely determined.

\medskip

Lemma 1 as a preliminary consideration connects (through Theorem D)
the Riemann hypothesis with the M\"obius function but it is not used in its proof.

\medskip

{\bf Lemma 1.} {\it Let $\varphi_n=\sum\limits^n_{k=1}h^{(n)}_ka_k$ be a sequence of functions in $\frak{R}$ that is convergent to $\chi$
with respect to the norm of $L_2^*$. Then the coefficient $h^{(n)}_k$ converges to $-\mu(k)$ as $n\to\infty$ for each $k\in\mathbb N.$}

\medskip

{\bf Proof.} Indeed, first by Remark 2 the function $\varphi_n$ is linear on each interval
$(l-1,l)$, $l=1,2,\ldots$, and its slope
\begin{equation}\label{present}
s_n\ =\ \sum\limits_{k=1}\limits^{n} k^{-1}h_k\ \to\ 0\ \ \
\hbox{as $n\to\infty$}
\end{equation}
because it should be
\begin{equation}\label{slope}
\lim\limits_{n\to\infty}\int\limits_0^1\varphi_n^2(t)\
\frac{dt}{t^2}\ =\ 0\ .
\end{equation}

Then we have that $h^{(n)}_1 \to -1=-\mu(1)$ as $n\to\infty$ because
$$
\lim\limits_{n\to\infty}\int\limits_1^2\left(1-\varphi_n(t)\right)^2\
\frac{dt}{t^2}\ =\
\lim\limits_{n\to\infty}\int\limits_1^2\left(1-\left(s_nt-h^{(n)}_1\right)\right)^2\
\frac{dt}{t^2}\ =\ 0\ ,
$$
again by Remark 2 and, similarly, $h^{(n)}_1+h^{(n)}_2\to0$ as $n\to\infty$ because
$$
\lim\limits_{n\to\infty}\int\limits_2^3\left(1-\varphi_n(t)\right)^2\,
\frac{dt}{t^2} =
\lim\limits_{n\to\infty}\int\limits_2^3\left(1-\left(s_nt-h^{(n)}_1-\left(h^{(n)}_1+h^{(n)}_2\right)\right)\right)^2\
\frac{dt}{t^2} = 0
$$ and, arguing by induction, we have that
$\sum\limits_{d|\,k}h^{(n)}_{d} \to 0$ for all $k\in\mathbb N$ because
$$
\lim\limits_{n\to\infty}\,\int\limits_k^{k+1}\,\left(\,1\, -\,
\varphi_n(t)\,\right)^2\ \frac{dt}{t^2}\ =
$$
$$
=
\lim\limits_{n\to\infty}\int\limits_k^{k+1}\left(1-\left(\left(\left(s_nt-h^{(n)}_1\right)-\ldots
-\sum\limits_{d_*\,|\,(k-1)}h^{(n)}_{d_*}\right)-\sum\limits_{d|\,k}h^{(n)}_{d}\right)\right)^2\
\frac{dt}{t^2} = 0\ .
$$
Thus, by Remark 3 we have that the coefficient $h^{(n)}_k$ converges to $-\mu(k)$ as $n\to\infty$ for each $k\in\mathbb N.$ $\Box$

\medskip

In the connection with the M\"obius function, let us recall also one more result of Baez-Duarte, see Theorem
2.2 in \cite{Baez2}, that after the replacement $t \mapsto 1/t$ can
be formulated in the following way.

\medskip

{\bf Theorem D$^*$.} {\it The sequence of the functions
\begin{equation}\label{eqBD}
\varphi_n\
:=\ -\sum\limits^{n}_{k=1}\mu(k)\,a_k
\end{equation}
cannot be convergent to $\chi$ with respect to the norm in $L^*_2$.}

\medskip

However, later on we effectively use only the following well--known identity
\begin{equation}\label{eq14}\boxed{\ \sum\limits_{n\leq t}\ \mu(n)\left[\frac{t}{n}\right]\ =\ 1 \ }\ \ \ \ \ \forall\
t\geq 1\end{equation}
see \cite{L}, p. 576, see also \cite{Ap}, p. 54, Example 1, \cite{HW}, Sec. 18.4, \cite{Te}, Sec. I.2.

\cc
\section{The main result}

{\bf Theorem 1.} {\it The Riemann hypothesis is true, i.e.,
the $\zeta -$function has as its zeros only complex
numbers with real part $1/2$ in addition to its trivial zeros at the
negative even integers.}

\medskip

{\bf Proof.} Let us consider the sequence of the functions in $\frak R$ of the form
\begin{equation}\label{eqfstar}
\varphi^*_n\
:=\ -\sum\limits^{n!-1}_{k=1}\mu(k)\,a_k\ +\ \delta(n)n!a_{n!}\ ,
\end{equation}
where
\begin{equation}\label{eqmstar}
\delta(n):=\sum\limits^{n!-1}_{k=1}\frac{\mu(k)}{k}\ ,
\end{equation}
and show that $\varphi^*_n$ converges by the norm in $L_2^*$ to the characteristic function $\chi$ of the interval $[1,\infty)$.

\medskip

Indeed, by Remark 2 the functions $\varphi_n^*\in\frak{R}$ are constant
on all open intervals of length $1$ that appeared in $(0,\infty)$ after removing all natural
numbers, continuous from the right at each point $k\in\mathbb N$, $\varphi_n^*\equiv 0$ on $(0,1)$, and by (\ref{eq14})
$\varphi_n^*\equiv 1$ on $[1,n!)$. Note also that functions $\varphi_n^*$ are $n!-$periodic because each function $a_k$, $k=1,\ldots , n!$ is $n!-$periodic. Thus, $\varphi_n^*\equiv 0$ on each interval $[n!k,n!k+1)$, $k\in \Bbb N$, and coincide with the function $\chi$ outside of all these intervals. Thus, the deviation of $\varphi_n^*$ from $\chi$ by the norm in $L_2^*$ is easy estimated
\begin{equation}\label{eqestar}
\| \chi - \varphi_n^*\|^2_* = \sum\limits^{\infty}_{k=1}\frac{1}{n!k(n!k+1)} < \frac{1}{n!^2}\sum\limits^{\infty}_{k=1}\frac{1}{k^2} = \frac{\zeta(2)}{n!^{2}} \to 0\ \ \ \hbox{as $n\to\infty$}\ .
\end{equation}

\medskip

Consequently, the conclusion of Theorem 1 follows from Theorem D. $\Box $

\bigskip

It should be emphasized here that the main result actually follows directly from the long known identity (\ref{eq14}) and the sufficient part of Theorem A. But Theorems D and D$^*$, as well as the geometric interpretation in Remark 2, and Lemma 1, serve merely as guiding considerations, i.e., as some clues for finding possible approximating sequences. Thus, in fact mathematicians already had everything they needed for this in 1950. There was just one last step to take.

\bigskip

{\bf Acknowledgments.} This research of the author during the period from 2024 to 2026 was partially supported
by the Simons Foundation, the grants PD-Ukraine-00010584 and SFI-PDUkraine-00017674.

\medskip
\noindent
{\bf Vladimir Ryazanov,} https://orcid.org/0000-0002-4503-4939,\\
Institute of Applied Mathematics and Mechanics,\\
National Academy of Sciences of Ukraine,\\
Institute of Mathematics, UKRAINE,\\
vl.ryazanov1@gmail.com,\\
Ryazanov@nas.gov.ua


\begin{thebibliography}{100}
\small


\bibitem{Ap} {\sc Apostol T.M.}, {\sl Introduction to analytic number
theory}, Undergraduate Texts in Mathematics, Springer-Verlag, New
York etc., 1976.


\bibitem{Ay} {\sc Ayoub R.}, {\sl An Introduction to the Analytic Theory of
Numbers}, Mathematical surveys {\bf 10}, AMS, Providence, 1963.

\bibitem{Baez1} {\sc Baez-Duarte L.}, {\sl On Beurling's real variable reformulation of the
Riemann hypothesis}, Adv. Math. {\bf 101} (1993), no. 1, 10--30.

\bibitem{Baez2} {\sc Baez-Duarte L.}, {\sl A class of invariant unitary operators}, Adv.
Math. {\bf 144} (1999), no. 1, 1--12.

\bibitem{Baez3} {\sc Baez-Duarte L.}, {\sl New versions of the Nyman-Beurling criterion for the Riemann
hy\-po\-the\-sis}, Int. J. Math. Math. Sci. {\bf 31} (2002), no. 7,
387–406.

\bibitem{Baez4} {\sc Baez-Duarte L.}, {\sl A strengthening of the Nyman-Beurling criterion for the Riemann
hypothesis}, Atti Accad. Naz. Lincei Cl. Sci. Fis. Mat. Natur. Rend.
Lincei (9) Mat. Appl. {\bf 14} (2003), no. 1, 5--11.

\bibitem{Baez5} {\sc Baez-Duarte L.}, {\sl A general strong Nyman-Beurling criterion for the Riemann
hypothesis}, Publ. Inst. Math. (Beograd) (N.S.) {\bf 78(92)} (2005),
117--125.

\bibitem{BBLS} {\sc Baez-Duarte L., Balazard M., Landreau B., Saias E.}, {\sl Notes sur la
fonction $\zeta$ de Riemann, 3}, (French. English summary) {\sl
Notes on the Riemann $\zeta$-function, 3}, Adv. Math. {\bf 149}
(2000), no. 1, 130–144.


\bibitem{Bag} {\sc Bagchi B.}, {\sl On Nyman, Beurling and Baez-Duarte's
Hilbert space reformulation of the Riemann hypothesis}, Proc. Indian
Acad. Sci. Math. Sci. {\bf 116} (2006), no. 2, 137--146.

\bibitem{BSl} {\sc Balazard M., Saias E.}, {\sl Notes sur
la fonction $\zeta$ de Riemann, 1}, (French) [{\sl Notes on the
Riemann $\zeta -$function. 1}], Adv. Math. {\bf 139} (1998), no. 2,
310--321.

\bibitem{BS2} {\sc Balazard M., Saias E.}, {\sl The
Nyman-Beurling equivalent form for the Riemann hypothesis}, Expo.
Math. {\bf 18} (2000), no. 2, 131--138.

\bibitem{Ber} {\sc Bercovici H., Foias C.}, {\sl A real variable restatement of Riemann's hypothesis}, Israel
J. Math. {\bf 48} (1984), no. 1, 57--68.

\bibitem{Beu} {\sc Beurling A.}, {\sl A closure problem related to the Riemann
zeta-function}, Proc. Nat. Acad. Sci. U.S.A. {\bf 41} (1955),
312--314.

\bibitem{B1} {\sc Broughan K.}, {\sl Equivalents of the Riemann hypothesis}, V.
1, Encyclopedia of Ma\-the\-ma\-tics and its Applications {\bf 164},
Cambridge University Press, Cambridge, 2017.

\bibitem{B2} {\sc Broughan K.}, {\sl Equivalents of the Riemann hypothesis}, V.
2, Encyclopedia of Ma\-the\-ma\-tics and its Applications {\bf 165},
Cambridge University Press, Cambridge, 2017.

\bibitem{B3} {\sc Broughan K.}, {\sl
Equivalents of the Riemann hypothesis. Vol. 3. Further steps towards
resolving the Riemann hypothesis Broughan}, Encyclopedia Math. Appl.
{\bf 187}, Cambridge University Press, Cambridge, 2024.

\bibitem{Bur} {\sc Burnol J.-F.}, {\sl On an analytic estimate in the theory of the Riemann
zeta function and a theorem of Baez-Duarte}, Acta Cient. Venezolana
{\bf 54} (2003), no. 3, 210--215.

\bibitem{Ch} {\sc Chen Ch.P.}, {\sl The Riemann
hypothesis and gamma conditions}, J. Math. Anal. Appl. {\bf 173}
(1993), no. 1, 258--275.


\bibitem{E} {\sc Edwards H.M.}, {\sl Riemann's zeta function}, Pure and
Applied Mathematics {\bf 58}, Academic Press, New York-London, 1974.


\bibitem{H} {\sc Habsieger L.} {\sl On the Nyman-Beurling criterion for
the Riemann hypothesis}, Funct. Approx. Comment. Math. {\bf 37}
(2007), part 1, 187--201.


\bibitem{HW} {\sc Hardy G.H., Wright E.M.} {\sl An Introduction to the Theory of
Numbers}, Fourth edition, Clarendon Press, Oxford, 1975.

\bibitem{Iv} {\sc Ivic A.}, {\sl The Riemann zeta-function. The theory of the Riemann
zeta-function with applications}, Wiley-Interscience Publication,
John Wiley Sons, Inc., New York, 1985.

\bibitem{Iw} {\sc Iwaniec H.}, {\sl Lectures on the Riemann zeta function}, University
Lecture Series {\bf 62}, American Mathematical Society, Providence,
RI, 2014.

\bibitem{KV} {\sc Karatsuba A.A.,
Voronin S.M.}, {\sl The Riemann zeta-function}, De Gruyter
Expositions in Mathematics {\bf 5}, Walter de Gruyter Co., Berlin,
1992.

\bibitem{L} {\sc Landau E.}, {\sl Handbuch der Lehre von der Verteilung der Primzahlen},
Teubner, Leipzig, 1909; third edition Chelsea, New York, 1974.

\bibitem{LR} {\sc Landreau B., Richard F.},  {\sl Le
critere de Beurling et Nyman pour l'hypothese de Riemann: aspects
numeriques}, [{\sl The Beurling-Nyman criterion for the Riemann
hypothesis: numerical aspects}], Experiment. Math. {\bf 11} (2002),
no. 3, 349--360.

\bibitem{Lap} {\sc Lapidus M.L.} {\sl In search of the Riemann zeros. Strings, fractal
membranes and noncommutative spacetimes}, American Mathematical
Society, Providence, RI, 2008.

\bibitem{Lee} {\sc Lee J.}, {\sl Convergence and the Riemann
hypothesis}, Commun. Korean Math. Soc. {\bf 11} (1996), no. 1,
57--62.

\bibitem{MR1} {\sc Maier H., Rassias M.Th.}, {\sl On the size of an expression in
the Nyman-Beurling-Baez-Duarte criterion for the Riemann
Hypothesis}, Canad. Math. Bull. {\bf 61} (2018), no. 3, 622--627.

\bibitem{MR2} {\sc Maier H., Rassias M.Th.}, {\sl Estimates of sums related to the Nyman-Beurling criterion for the
Riemann hypothesis}, J. Number Theory 188 (2018), 96--120.

\bibitem{Nik} {\sc Nikolski N.}, {\sl Distance
formulae and invariant subspaces, with an application to
localization of zeros of the Riemann $\zeta$-function}, Ann. Inst.
Fourier (Grenoble) {\bf 45} (1995), no. 1, 143–159.

\bibitem{Ny} {\sc Nyman B.}, {\sl On the One-Dimensional Translation Group and
Semi-Group in Certain Function Spaces}, Thesis, University of
Uppsala, 1950. 55 pp.

\bibitem{R} {\sc Riemann B.}, {\sl \"Uber die Anzahl der Primzahlen unter einer gegebenen
Gr\"osse}, Monat. der Konigl. Preuss. Akad. der Wissen. zu Berlin
aus der Jahre 1859 (1860), 671–680; see also Gesammelte Werke
(Teubner, Leipzig, 1892; Dover, New York, 1953); see the English
transl. in http://
www.maths.tcd.ie/pub/HistMath/People/Riemann/Zeta/EZeta.pdf

\bibitem{Ro} {\sc Roton A.}, {\sl Generalization du critere de
Beurling-Nyman pour l'hypothese de Riemann} [{\sl Generalization of
the Beurling-Nyman criterion for the Riemann hypothesis}], Trans.
Amer. Math. Soc. {\bf 359} (2007), no. 12, 6111--6126.

\bibitem{Te} {\sc Tenenbaum G.}, {\sl Introduction to analytic and probabilistic number theory}, 3rd ed.,
Grad. Stud. Math., {\bf 163}, American Mathematical Society, Providence, RI,  2015.

\bibitem{T} {\sc Titchmarsh E.C.}, {\sl The Theory of the Riemann Zeta-Function},
Clarendon Press, Oxford, 1951.

\bibitem{Ya2} {\sc Yang J.}, {\sl A generalization of Beurling's criterion for the Riemann
hypothesis}, J. Number Theory {\bf 164} (2016), 299--302.

\bibitem{Ya} {\sc Yang J.}, {\sl A Friedrichs angle between the Nyman-Beurling spaces and the Riemann
hypo-thesis}, J. Math. Anal. Appl. {\bf 560} (2026), no. 1, Paper No. 130494, 13 pp.

\bibitem{Ya3} {\sc Yang J.}, {\sl A generalization of Weingartner's theorem related to the Riemann hypothesis},
Arch. Math. (Basel) {\bf 126} (2026), no. 6, 625–633.

\end{thebibliography}
\end{document}